\documentclass[reqno,11pt,a4]{amsart}
\usepackage{amsmath,amssymb,tabularx,setspace,color,graphicx}
\makeatletter
\renewcommand\section{\@startsection {section}{1}{\z@}%
                                   {-3.5ex \@plus -1ex \@minus -.2ex}%
                                   {2.3ex \@plus.2ex}%
                                   {\normalfont\large\bfseries}}
\makeatother
\makeatother

\newtheorem{Theorem}{Theorem}
\usepackage{amsmath}
\usepackage{lipsum}
\usepackage[margin=3.5cm]{geometry}
\usepackage{amsmath}
\usepackage{setspace}
\usepackage{amsthm}
\usepackage{graphicx}
\usepackage{enumerate}
\usepackage{multicol}
\usepackage{graphics}
\usepackage{bigstrut}
\usepackage{multirow}
\usepackage{bigstrut}
\usepackage{tabularx}
\usepackage{setspace}
\usepackage{amssymb}
\usepackage{amsmath,amssymb,color}
\usepackage[figuresright]{rotating}
\usepackage{enumerate}
\usepackage{graphics}
\usepackage{graphicx}
\usepackage{multirow}
\usepackage{lscape}
\usepackage[dvipsnames]{xcolor}
\usepackage{mathrsfs}
\usepackage{hyperref}
\usepackage{textcomp}
\usepackage{listings}
\makeatletter

\begin{document}
		\title[LT models for competing risks data with cure fraction]{Semiparametric transformation models for competing risks data with cure fraction}
			\author[S\lowercase{udheesh} K. K.,  S\lowercase{reedevi} E. P. \lowercase{and} S\lowercase{ankaran} P. G.]{S\lowercase{udheesh} K. K\lowercase{attumannil}$^{\lowercase{a},\dag}$, S\lowercase{reedevi} E. P.$^{\lowercase{b}}$ \lowercase{and }S\lowercase{ankaran} P. G.$^{\lowercase{c}}$ \\
		$^{\lowercase{a}}$I\lowercase{ndian} S\lowercase{tatistical} I\lowercase{nstitute}, C\lowercase{hennai}, I\lowercase{ndia},\\
			$^{\lowercase{b}}$SNGS C\lowercase{ollege}, P\lowercase{attambi}, I\lowercase{ndia},\\
$^{\lowercase{c}}$C\lowercase{ochin} U\lowercase{niversity of} S\lowercase{cience and} T\lowercase{echnology}, K\lowercase{ochi}, I\lowercase{ndia.} \\
}
\thanks{${\dag}$ Corresponding author E-mail: skkattu@isichennai.res.in}
\doublespace
    \begin{abstract} We propose a new method for the analysis of competing risks data with long term survivors. The proposed method enables us to estimate the overall survival probability and cure fraction simultaneously. We formulate the effect of covariates on cumulative incidence functions using linear transformation models. Estimating equations based on counting process are developed to estimate regression coefficients. The asymptotic properties of the estimators are studied using martingale theory. An extensive Monte Carlo simulation study is carried out to assess the finite sample performance of the proposed estimators.  Finally, we illustrate our method using a real data set. 
\end{abstract}
\maketitle
\vspace{-0.2in}
\section{Introduction}
\vspace{-0.1in}
\noindent Competing risks emerge naturally in lifetime data analysis when the subjects under study are at risks for more than one cause of failure. For example, consider a study of patients suffering from heart  disease. The patients may die due to other causes like an  accident or other diseases, which may alter the probability of death due to heart disease. The cause specific hazard  and cumulative incidence  functions are commonly employed for the analysis of competing risks data. For a comprehensive review on this topic, one may refer to Crowder (1997), Kalbfleisch and Prentice (2002) and Lawless (2002).

In lifetime studies, we generally assume that all subjects under study will experience the event of interest, if they are followed sufficiently long. However, in some situations, a non-negligible proportion of individuals may not experience the event of interest  even after a long period of observation. For example, in many clinical trials, there exist a proportion of subjects who may not experience the event of interest. These patients can not be treated  as censored, and should be considered as cured individuals. In such a scenario, traditional methods for analyzing survival data have to be modified. Cure rate models have wide range of applications in many fields including medical and public health. In cure rate models, the entire population is considered as a mixture of two groups of patients, susceptible (non-cured) and non-susceptible (cured). In clinical trials of cancer studies, the proportion of cured individuals is an important factor in estimating survival probabilities. Hence, it is of interest to develop models which incorporate cure fraction.

The mixture cure rate model (MCM) proposed by Boag (1949) has been popular in the analysis of survival data with long term survivors. Some important works in this area include Farewell (1982), Taylor (1995), Peng and Dear (2000) Sy and Taylor (2000) and  Zhang and Peng (2007) among others.  The books by Maller and Zhou (1996) and Ibrahim, Chen  and Sinha (2002)  served as  excellent references on this topic.

Competing risks data with cure fraction arises when the susceptible group is exposed to different risks of failure. For example, consider the data obtained from a clinical trial on HIV infection and AIDS of 329 homosexual men from Amsterdam, reported in Geskus (2015).  During the course of period from HIV infection (NSI phenotype) to death, the virus may be switched to SI phenotype or AIDS; which are mutually exclusive events and both affect the subsequent disease progression differently. Accordingly it is important to observe which infection has occurred first. The presence or absence of the CCR5-$\Delta$32 deletion in one or both chromosomes is considered as the covariate in the study. An interesting characteristic of the data is that time of occurrence of these events vary from 0.112 years (approximately 41days) to 13.936 years (approximately 5090 days).  We also found that after 12.4 years (4526 days), there are only 3 events observed (one due to AIDS and two due to SI), whereas around 14.5\% of the total lifetimes reported are larger than 12.4 years. This scenario indicates the presence of cured individuals in the data (Maller and Zhao, 1996). Hence, it is required to analyse the data using cure rate model with competing risks. In Section 5, we use this data to illustrate the applicability of the proposed method.

Modelling and analysis of competing risks data with cure fraction is considered by several authors in literature. A Bayesian method that  unifies the mixture cure and competing risks approach is developed by Basu and Tiwari (2010). Choi and Huang (2014) considered a finite mixture model for analysing competing risks with cure fraction.  Choi, Huang, and Cromier (2015) proposed a  semiparametric mixture model to analyse competing risks data with cure fraction using a multinomial logistic model.  A semiparametric accelerated failure time model for the cause-specific survival function that is combined through a multinomial logistic model within the cure-mixture modeling framework is developed by Choi, Zhu and Huang (2018).  The vertical modeling approach is extended to analyse competing risks data with cure fraction by Nicolaie, Taylor and Legrand (2019). Recently,  Wang, Zhang and Tang (2020) have proposed a semiparametric estimation procedure for the  accelerated failure time mixture cure model in the presence of competing risks. However, in all the aforementioned and other related works, one needs to estimate the probability of failure due to each cause in the presence of covariates, in order to estimate the overall survival function.

Motivated by this, we propose a new semi-parametric regression model to analyse competing risks data with cured individuals. We use a mixture model approach to analyse the competing risks data with cure fraction. We specify cumulative incidence  functions of the competing risks data using linear transformation model.  Unlike the traditional approach discussed earlier, the novelty of our approach is that we can estimate the overall survival function without estimating the cure fraction separately. This achievement is in the line of promotion  time cure model discussed by  Tsodikov (2002) in the non-competing risks scenario.

The rest of the  article is organized as follows. In Section 2, we propose a semiparametric regression model for the analysis of  competing risks data with cure fraction.  In Section 3, we develop counting process based estimating equations to find the estimators of regression coefficients. The asymptotic properties of the estimators are also studied.   An extensive Monte Carlo simulation study is carried out to assess the finite sample performance of the proposed  estimators.  Computational algorithms along with the results of the simulation study are presented in Section 4. In Section 5, we illustrate the application  of the proposed method using a real data set obtained from a study of HIV and AIDS infection. Finally,  in Section 6, we summarise  major conclusions of the study along with discussion on some open problems.
\section{Model and methods}
We propose a new method to analyse  competing risks data with cure fraction. Consider a general competing risks setup with $K$ distinct failure types. Let $(T,J)$ be the observed  data,   where $T$ denotes the time to failure and  $J\in \{0, 1, … , K\}$ be the corresponding cause of failure. We use  $J = 0$ to denote a  subject who is insusceptible to any type of failures. The mixture modeling approach has been popular in the analysis of lifetime data with long term survivors. In the presence of cure fraction, mixture model of competing risks data assumes that the failure time $T$ can be decomposed as
\[T = \sum_{k=1}^{K}T_k .I(J = k)+ \infty. I(J = 0),\]
where $T_k$ denote the  latent failure times due to cause $k$, $k=1,\ldots,K$ and  $I$  denotes the indicator function.

Let $Z$  be $p$-dimensional covariates, possibly time variant.  We assume that the censoring random variable $C$ is  independent of $T$, conditional on the covariates $Z$. 
Denote $\pi_k(Z)=P(J=k|Z)$ as the probability  of experiencing the event from  cause $k,\, k=1,\ldots,K$ and $\pi_0(Z)$ as the cure fraction. Clearly, $\pi_0(Z)=1-\sum_{k=1}^{K}\pi_k(Z)$.
Mixture cure rate model assumes that the overall conditional survival function
of $T$ given $Z$ has the  representation given by
\begin{equation}\label{mixsur1}
S(t|Z) = P(T > t | Z) =\pi_0(Z)+\sum_{k=1}^{K}\pi_k(Z)S_k(t|Z),
\end{equation}where $S_k(t|Z)=P(T_k>t|J=k,Z)$ and $\beta=\{\beta_1\ldots \beta_k\}$ denote the full set of parameters.  To estimate $S(t|Z)$ specified in (\ref{mixsur1}), we need to model the  effect of covariates on $\pi_k(Z)$ and $S_k(t|Z)$ separately (Patilea and Keilegom, 2017).  Under the finite-mixture model, apart from estimating $S_k(t|Z)$, practitioners usually model $\pi_k(Z)$ using logistic regression (Farewell, 1982).  This makes the implementation of the method complex and computationally challenging. Motivated by this, we propose a new method to analyse the competing risks data with long term survivors.

We consider a semi-parametric transformation  model for cumulative incidence functions, which is the probability of failure from any one cause in the presence of other risks over a certain time period. The cumulative incidence function conditional on $Z$ denoted by  $F_k(t|Z)$ is given by
$$F_k(t|Z)=P(T\le t,J=k|Z),\,\,k=1,\ldots,K.$$
It is the cumulative probability of failure from the $j$-th cause of failure in the presence
of the remaining causes of failure, conditional on the covariates.
Now, we express $S(t|Z)$ given in (\ref{mixsur1})  in terms of $F_k(t|Z)$. For $k=1,2,...,K$, consider
\begin{eqnarray}\label{conditional}
	F_k(t|Z)&=&P(T\le t,J=k|Z)\nonumber \\	&=&P(T\le t|J=k,Z)P(J=k|Z)\nonumber\\&=&(1-P(T>t|J=k,Z))P(J=k|Z)\nonumber\\&=& \pi_k(Z)(1-S_k(t|Z)).
\end{eqnarray}
The representation given in (\ref{conditional}) allows us to model the conditional cumulative incidence function through $\pi_k(Z)$. We use linear transformation models to specify  each cumulative incidence function $F_k(t|Z)$ without using $\pi_k(Z)$. Thus, using (\ref{conditional}), we rewrite $S(t|Z)$  given in (\ref{mixsur1}) as
\begin{eqnarray}
\label{mixsurnew}
S(t|Z)
&=&1-\sum_{k=1}^{K}F_k(t|Z).
\end{eqnarray}


In the proposed method, we  formulate the conditional cumulative function $F_k(t|Z)$ specified in (\ref{mixsurnew}), through a class of linear transformation models proposed by (Mao and Lin, 2017)
\begin{equation}\label{LTMCOMP}
g_k\{F_k(t|Z)\} = h_k(t) + Z'\beta_k,\quad k=1,2,\ldots,K,
\end{equation}where $g_k$ is  a known increasing cause specific link function, $h_k(t)$  is an unknown non-decreasing function and $\beta_k$ is a set of p-dimensional regression coefficients. We assume $h_k(0)=-\infty$.  The  proportional hazards (PH) model and proportional odds (PO) model are special cases of (\ref{LTMCOMP}) when $g_k(x)$ takes the form $\log(-\log(1-x))$ and $\log(x/(1-x))$, respectively. For more details on linear transformation model one can refer to Doksum (1987), Chen, Jin and Ying (2002), Zeng and Lin (2006) and Mao and Lin (2017). The LT model in (\ref{LTMCOMP}) can be represented  as
\begin{equation}\label{LTnew}
  h_k(T)= -Z'\beta_k+\epsilon_k,\quad k=1,2,\ldots,K,
\end{equation}
where $\epsilon_k$ is a random error with
a known distribution $F_{\epsilon_{k}}(x)=P(\epsilon_k \le x)=g_k^{-1}(x)$ and is independent of the covariate  $Z$. The model in (\ref{LTnew}) reduces to PH and PO model, when $\epsilon_k$  has standard extreme value and standared logistic distribution, respectively.   For further details of this equivalent form,  see  Fine, Ying and Wei (1998), Zhang,  Sun, Zhao and Sun  (2005) and Guo and Zeng (2014).
Define conditional cause specific hazard rate  function as
$$\lambda_k(t|Z)=\lim_{\delta t\rightarrow 0}\frac{1}{\delta t} Pr\left(t\le T<t+\delta t,J=k|T\ge t,Z\right).$$
Clearly, $F_k(t|Z)=1-\exp(-\Lambda_k(t|Z)),$  where $ \Lambda_k(t|Z)= \int_{0}^{t}\lambda_k(s|Z)ds$.
    We use the form given in (\ref{LTnew}), while developing the estimating equations based on Martingale. Hence,  inference on
$F_k(t|Z)$  follows from the simple relationship $F_k(t|Z)=1-\exp\left(-\Lambda_{k}(t|Z)\right)$.

 From  (\ref{LTMCOMP}), we have
\begin{equation}\label{LTMCOMPSUR}
  F_k(t|Z) =g_k^{-1}(h_k(t) + Z'\beta_k),\quad k=1,2,\ldots,K.
\end{equation}
Thus using(\ref{mixsurnew}), the survival function of $T$ given $Z$  is given by
\begin{equation}\label{ovsurv}
S(t|Z)=1-\sum_{k=1}^{K}g_k^{-1}(h_k(t) + Z'\beta_k).
\end{equation}
In view of (\ref{cureprob}) given in Section 3, we can easily see that as $t\rightarrow\infty$, $S(t|Z)$ become $\pi_0(Z)$, the cure fraction. From the representation (\ref{ovsurv}), it is obvious that we can find the overall survival function of $T$  without estimating the cure fraction.  For finding overall survival function $S(t|Z)$, we need to estimate $h_k(.)$ and $\beta_k$ in (\ref{ovsurv}). The estimation  of $h_k(t)$ and $\beta_k$ are discussed in the next section.


\section{Estimating equation and asymptotic properties}
\noindent
Let $T$ and $C$ be the failure time and the censoring time random variables, respectively. Let $\tilde T=\min(T,C)$, then we observe $\tilde T$ and $\delta$ where $\delta=I(T<C)$ with $I$ as the indicator function. Let $J$ be a discrete random variable with support $\{0,1,\ldots,K\}$. We assume that the failure of a subject  is due to any of the causes $\{1,2,\ldots,K\}$ and we
use the notation $J = 0$ to denote the subject which are insusceptible to any of the $K$ causes. Let $Z$ be $p$-dimensional covariate vector.  The observed data $(\tilde T_i,\delta_i,\delta_i J_i, Z_i)$, $i=1,2,\ldots,n$, is the independent
copies of the vector $(\tilde T,\delta,\delta J, Z)$.  Let $\lambda_k(t)$ and $\Lambda_k(t)$, $k=1,\ldots,K$  be the cause specific hazard function and  cause specific cumulative  hazard function conditional on $Z$ corresponding to $k$-th cause.

To estimate $\beta_k$ and $h_k(.)$, $k=1,\ldots,K$ we propose the estimating equations based on  counting process. Let $N_{ik}(t)=\delta_i.I(\tilde T_i\le t,J=k)$  be the counting process associated with failure time of $i$-th subject due to cause $k$, $k=1,\ldots,K$ and $i=1,2,\ldots,n$. Let $N_k(t)=\sum_{i=1}^{n}N_{ik}(t)$ be the number of failures  due to cause $k$ by time $t$. Then $N(t)=\sum_{k=1}^{K}N_{k}(t)$ is the total number of events experienced by time $t$. Define  at-risk process $Y_i(t)=I(\tilde T_i>t)$ and $Y(t)=\sum_{i=1}^{n}Y_i(t)$.  It can be easily verified that $N_{ik}(t)$, $i=1,\ldots,n$ and $k=1,\ldots,K$
are  local sub-martingales with respect to appropriate filtration (Andersen et al., 1993). By Doob-Meyer decomposition, martingale process associated with $N_{ik}(t)$ is given by (Andersen et al., 1993)
$$M_{ik}(t)=N_{ik}(t)-\int_{0}^{t}Y_i(t)d\Lambda_{\epsilon_{k}}(t|Z),\,\,i=1,\ldots,n,\, k=1,2,\ldots,K,$$
where $\Lambda_{\epsilon_{k}}$ is the cumulative cause-specific hazard functions of $\epsilon_{k}$.
Using linear transformation model  specified in (\ref{LTnew}), we have
\begin{equation}\label{martingale}
M_{ik}(t)=N_{ik}(t)-\int_{0}^{t}Y_i(t)d\Lambda_{\epsilon_{k}}(h_k(t) + Z'\beta_k).
\end{equation}
By definition,  $M_{ik}(t)$ is a mean zero martingale process with respect to appropriate
filtration (Andersen and Gill, 1982). We now propose  the  following estimating equations to obtain the estimators of $\beta_k$ and $h_k(.)$, $k=1,2,\ldots,K$,
\begin{equation}\label{betahat}
U_{\beta_k}(\beta_k,h_k)=\sum_{i=1}^{n}\int_0^\infty {Z_i} \big[dN_{ik}(u)- Y_i(u)d\Lambda_{\epsilon_{k}}(h_k(t) + Z_i'\beta_k)\big]=0,
\end{equation}and
\begin{equation}\label{LTMO}
U_{h_k}(\beta_k,h_k) =\sum_{i=1}^{n}\big[dN_{ik}(t)- Y_i(t)d\Lambda_{\epsilon_{k}}(h_k(t) + Z'\beta_k)\big]=0,\,t\geq 0.
\end{equation}Solving the above equations iteratively   we obtain the estimators of $\beta_k$ and $h_k(.)$, $k=1,2,\ldots,K$. The estimating  equations in (\ref{betahat}) and (\ref{LTMO}) reduce to the estimating equations given by Chen et al. (2002) when $k=1$.

When the estimators of $\beta_k$ and $h_k(t)$ are obtained, one can estimate the overall survival function without estimating cure probability. However, our newly proposed method enables us to estimate the cured probability from the proposed model itself.  We use  (\ref{LTMCOMPSUR})  to estimate the cure fraction. From the  definition of cumulative incidence function we have (Lawless, 2002),
$$P(J=k|Z)=F_k(\infty,J=k|Z)=g_k^{-1}(h_k(\infty)+Z'\beta_k).$$ Hence, the cure probability can be obtained using the relation
\begin{equation}\label{cureprob}
  \pi_0(Z)=1-\sum_{k=1}^{K}g_k^{-1}( h_k(\infty)+Z'\beta_k).
\end{equation}
Thus the cure fraction is estimated  by
\begin{equation}\label{cureprobest}
  \widehat\pi_0(Z)=1-\sum_{k=1}^{K}g_k^{-1}(\widehat  h_k(\infty)+Z'\widehat\beta_k).
\end{equation}
In practice, $h_k(\infty)$ can be estimated by $\max(\widehat h_k(t))$. We now provide the asymptotic distribution of the $\widehat\beta_k$. To find the distribution of $\widehat\beta_k$, first we show that each $\widehat{h}_k(t)$ is a consistent estimator of $h_k(t)$. Then we use the estimating equation in (\ref{betahat}) to  find the asymptotic
distribution of $\widehat\beta_k$. We assume the following regularity conditions to prove the asymptotic properties.
\begin{enumerate}[D1.]
	\item[D1.] The covariates $Z$ are bounded in probability.
	\item[D2.]  The true value of the parameters $\beta_k$, $k=1,\ldots, K$ lies in a compact set of ${R}^p$
	\item[D3.] Define $\tau=\inf\{t:P(\tilde{T}>t)=0\}$. 
	\item[D4.]  The derivatives of $\lambda_{\epsilon_{k}}(.),\,k=1,\ldots,K$ exists and continuous.
\item[D5.]  The  martingales defined in equation (\ref{martingale}) satisfies the  regularity conditions as in Fleming and Harrington (1991).
\end{enumerate}
The conditions $D1-D4$ are standard regularity conditions used in survival analysis. The assumption $D5$ is used to establish martingale central limit theorem.

The following additional notations are needed for deriving asymptotic distribution.
For any $s<t \in (0, \tau]$, define
\begin{eqnarray*}
	\noindent B_k(t,s)&=&\exp\Big(\int_{s}^{t}\frac{E[\partial\lambda_{\epsilon_{k}}/\partial t (Z_i'\beta_{k}+h_{k}(u))Y(u)]}{E[\lambda_{\epsilon_{k}} (Z_i'\beta_{k}+h_{k}(u))Y(u)]}dh_{k}(u)\Big).
\end{eqnarray*}
For $k=1,\ldots,K$, define $ \mu_k(t)=\frac{C_{zk}(t)}{C_{dk}(t)}$ where  
$$ C_{zk}(t)=E[Z\lambda_{\epsilon_{k}} (Z'\beta_{k}+h_{k}(t))Y(t)B_k(t,T)]$$ and $$C_{dk}(t)=E[\lambda_{\epsilon_{k}} (Z'\beta_{k}+h_{k}(t))Y(t)].$$

The asymptotic properties of  $\widehat\beta_k$ and $\widehat{h}_k(t)$ are established in the following theorems.
\begin{Theorem}
  Under the regularity conditions $D_1-D_4$, for  $k=1,\ldots,K,$  $\widehat\beta_k$ and $h_k(t)$ are strongly  consistent. That is
  $$\|\widehat\beta_k-\beta_{k}\|+\sup_{t\in [0,\tau]}\sum_{k=1}^{K}|\widehat h_k(t)-h_{k}(t)|\rightarrow0,$$
  almost surely, where $\|.\|$ denotes the Euclidean norm.

\end{Theorem}
\begin{Theorem}
Under the regularity conditions $D_1-D_5$, for  $k=1,\ldots,K$, as $n\rightarrow\infty$,  $\sqrt{n} (\widehat\beta_k-\beta_k)$ converges in distribution to multivariate normal with    zero mean vector  and variance-covariance matrix $\Sigma_k$ where $\Sigma_k=\Sigma_{1k}^{-1}\Sigma_{0k}(\Sigma_{1k}^{-1})'$ with
	\begin{equation}
\label{var}
\Sigma_{0k}=\int_{0}^{\tau}E\{({Z}-\mu_k(t))({Z}-\mu_k(t))'Y(t)
d\Lambda_{\epsilon_{k}}[Z_i'\beta_{k}+h_{k}(t)]\}
	\end{equation}
\begin{equation}\label{limvar}		
\Sigma_{1k}=\int_{0}^{\tau}E\{({Z}-\mu_k(t))(   {Z'}\partial\lambda_{\epsilon_{k}}/\partial t\{Z_i'\beta_{k}+h_{k}(t)\})Y(t)\}dh_{k}(t).
\end{equation}
\end{Theorem}
Proofs of Theorem 1 and Theorem 2 are given in Appendix.

\section{Computational Algorithm and Simulations}
The estimators $\widehat\beta_k$  and $\widehat h_k(t)$ are obtained as the solutions of  the equations  (\ref{betahat}) and (\ref{LTMO}). The value of $\widehat h_k(t)$ are estimated at observed failure time due to cause $k$.  For computational simplicity we express the set of equations (\ref{LTMO}) in an alternative form.   Let $t_{k1},t_{k2},\ldots,t_{km}$ be
the observed failure times due to the cause $k$, $k=1,2,\ldots,K$. Then (\ref{LTMO}) can be rewritten as
\begin{equation}\label{eq5.1}
\sum_{i=1}^{n}Y_i(t_{k1})\Lambda_{\epsilon_{k}}(Z_i'\beta_k+h_k(t_{k1}))=1,
\end{equation}
\begin{equation}\label{eq5.2}
\sum_{i=1}^{n}Y_i(t_{kj})\big(\Lambda_{\epsilon_{k}}(h_k(t_{kj})+Z_i'\beta_k)-\Lambda_{\epsilon_{k}}(h(t_{kj}-)+Z_i'\beta_k)\big)=1,\, j=2,3,\ldots,m.
\end{equation}
Thus, we have the following  iterative algorithms for computing  $\widehat\beta_k$ and $\widehat{h}_k(t)$.
\begin{enumerate}[Step 1.]
\item[Step 1.]  Choose an initial value ${\beta}_{k}^{(0)}$ for $\beta_k$,  $k=1,\ldots,K$. Obtain an estimator $\widehat{h}_k(t)$ for $h_k(t)$ by solving the equations (\ref{eq5.1}) and (\ref{eq5.2}). 
\item[Step 2.]  Find $\widehat{\beta}_k$  by solving the equations (\ref{betahat}) using $\widehat{h}_k(t)$ obtained in Step 2.
\item[Step 3.]  Set $\beta_{k}^{(0)}=\widehat{\beta}_k$ (the estimator obtained in the previous step) and repeat the Steps 1-3 until the convergence of $\widehat\beta_k$.
\end{enumerate}
We conduct an extensive Monte Carlo simulation study to assess the finite sample performance of the proposed estimators. In simulation, we consider two models where the cause specific  hazard function is specified by $\lambda_{\epsilon_{k}}(t)=\frac{e^{t}}{1+re^{t}},\,r=0,1,\,k=1,2$ and $h_k(T)=\log T_k$. Here, $r=0,1$  corresponds to the proportional hazards and the proportional odds model, respectively.

 The covariate $Z_1$ is generated from Bernoulli distribution with probability of success equal to $0.5$.  When $r=0$, we generate   $T$ using the model expression $-\log(U)\exp(-Z_1b_k)$, where $U$ is $U(0,1)$ random variable.  Also, when $r=1$, we simulate $T$  using the expression $((1-U)/U)\exp(-Z_1b_k)$.  In both the cases,  we choose the values of the regression parameters as $(\beta_{1},\beta_{2})=(2,2),(2,3),(3,2)$ and $(3,3)$.  Censoring random variable $C$ is simulated from $U(0,c)$ distribution where $c$ is chosen so that the  sample contain desired percentage  of censored observations, ie. $P(T>C)=q$, $0\le q<1$. We consider two different censoring scenarios with $q=\,0.2,\,0.4$.
We simulate observations with different  samples  sizes $n=100,200,500$ and simulation is repeated ten thousand times. The  simulation  is carried out using R program.


\begin{table}[h]\label{ph-bias and mse}
\caption{ Absolute bias and MSE of the estimators of regression coefficients under PH model}
\label{t:one}
\begin{small}
\centering
\begin{tabular}{crrrrrrrrr}

	\hline {Censoring(\%)}  & $ n$
&Bias&MSE&Bias&MSE&Bias&MSE&Bias&MSE\\
		\hline
		&&\multicolumn{2}{c}{$\beta_1=2$}&\multicolumn{2}{c}{$\beta_2=2$}&
\multicolumn{2}{c}{$\beta_1=2$}&\multicolumn{2}{c}{$\beta_2=3$}\\	
	\hline	
	\multirow{3}{*}{20} & 100 & 0.1227&	0.0309&	0.1197	&0.0378&0.1223&	0.0314&	0.0405&	0.0357 	\\
	& 200 &0.1150&	0.0232	&0.1098&	0.0256&0.1152&	0.0231&	0.0322&	0.0181 \\
	& 500 & 0.1016&	0.0189&	0.0922&	0.0197&0.1045&	0.0189&	0.0131	&0.0086
	\\
	\hline
	\multirow{3}{*}{40} & 100& 0.1226&	0.0317	&0.1185&	0.0465&0.1223&	0.0320&	0.0383&	0.0475 \\
	& 200 & 0.1150&	0.0237&	0.1039&	0.0297&0.1126&	0.0234	&0.0271&	0.0235
	\\
	& 500 &0.1011&	0.0192&	0.0834&	0.0217&0.1017&	0.0191&	0.0012&	0.0104
	\\
		\hline
		&&\multicolumn{2}{c}{$\beta_1=3$}&\multicolumn{2}{c}{$\beta_2=2$}&
\multicolumn{2}{c}{$\beta_1=3$}&\multicolumn{2}{c}{$\beta_2=3$}\\
		\hline
\multirow{3}{*}{20} & 100 &0.0428&	0.0256	&0.1212	&0.0368&0.0431&	0.0257&	0.0418	&0.0358\\
& 200 &0.0360&	0.0139&	0.1097&	0.0261&0.0363&	0.0139&	0.0303&	0.0182 \\
& 500 &0.0198&	0.0072	&0.0791&	0.0212&0.0226&	0.0067&	0.0112&	0.0085 \\
\hline
\multirow{3}{*}{40} & 100& 0.0431&	0.0278&	0.1160&	0.0486&0.0419&	0.0282&	0.0391&	0.0498\\
& 200 &0.0363&	0.0139&	0.0303&	0.0182&0.0363&	0.0150	&0.0242&0.0240
\\
& 500 &0.0198&	0.0072&	0.0791&	0.0212&0.0222&	0.0072&	0.0004	&0.0105
\\
\hline
\end{tabular}
\end{small}
\end{table}%

\begin{table}[h]\label{po-bias and mse}
\caption{ Absolute bias and MSE of the estimators of regression coefficients under PO model}
\label{t:one}
\begin{small}
\centering
\begin{tabular}{crrrrrrrrr}

	\hline {Censoring(\%)}  & $ n$
&Bias&MSE&Bias&MSE&Bias&MSE&Bias&MSE\\
		\hline
		&&\multicolumn{2}{c}{$\beta_1=2$}&\multicolumn{2}{c}{$\beta_2=2$}&
\multicolumn{2}{c}{$\beta_1=2$}&\multicolumn{2}{c}{$\beta_2=3$}\\	
	\hline	
	\multirow{3}{*}{20} & 100 & 0.1223	&0.0971	&0.0307&	0.0370 &0.1218&	0.0408	&0.0310&	0.0358 	\\
	& 200 &0.1161&	0.0931&	0.0235&	0.0259 &0.1149	&0.0303&0.0229&	0.0184  \\
	& 500 & 0.1033&	0.0208&	0.0189&	0.0200 &0.1042&	0.0167&	0.0188	&0.0084 \\
	\hline
	\multirow{3}{*}{40} & 100& 0.1213&	0.1183&	0.0317&	0.0476 &0.1218	&0.0382&	0.0319&	0.0487  \\
	& 200 & 0.1141&	0.1047&	0.0238&	0.0305 &0.1148	&0.0235	&0.0241&0.0231 	\\
	& 500 & 0.1004&	0.0832	&0.0189	&0.0215& 0.1003	&0.0050&0.0190	&0.0102 \\

		\hline
		&&\multicolumn{2}{c}{$\beta_1=3$}&\multicolumn{2}{c}{$\beta_2=2$}&
\multicolumn{2}{c}{$\beta_1=3$}&\multicolumn{2}{c}{$\beta_2=3$}\\
		\hline
\multirow{3}{*}{20} & 100 & 0.0440&	0.1193&	0.0258&	0.0371 &0.0427&	0.0419&	0.0265&	0.0359 \\
& 200 &0.0334&	0.1102&	0.0136&	0.0261 &0.0368&	0.0310&	0.0137&	0.0181   \\
& 500 &0.0234&	0.0936&	0.0068&	0.0197& 0.0229&	0.0138&	0.0069&	0.0084   \\
\hline
\multirow{3}{*}{40} & 100& 0.0411&	0.1183&	0.0275&	0.0471 &0.0430&	0.0387&	0.0276&	0.0491   \\
& 200 & 0.0346&	0.1027&	0.0145&0.0303 &0.0363&0.0256&0.0145&	0.0235 \\
& 500 & 0.0197&	0.0763&	0.0070&	0.0215 &0.0204&	0.0010&	0.0071&	0.0105 \\
\hline
\end{tabular}
\end{small}
\end{table}%

The regression parameters are estimated using the iterative algorithm given in the beginning of this section. The absolute bias and MSE (Mean Square Error) of the estimators of regression coefficients obtained under PH and PO models for different parameter settings  are reported in Tables 1 and 2. From Tables 1 and 2, we observe that the absolute bias of the estimators of regression coefficients  approaches  zero as sample size increases. In all cases, the MSE also decreases as sample size increases. We also note that the absolute bias and MSE of the estimators increase with the censoring percentage.

The absolute bias and the MSE of the estimates of $h_j(\cdot)$ for $j=1,2$ are estimated at three different time points $t=1,2,3$. The results are reported in Tables 3-6. In Tables 3 and 4, we present absolute bias and MSE of the estimators of $h_j(\cdot)$ for $j=1,2$ and for parameter combinations $(\beta_1,\beta_2)=(2,2)$ and $(3,2)$ under PH model. Under PO model, absolute bias and MSE of the estimators of $h_j(\cdot)$ for $j=1,2$  and for parameter combinations $(\beta_1,\beta_2)=(2,3)$ and $(3,3)$ are given in Tables 5 and 6. We obtain similar results for other parameter combinations also under both PH and PO model and we report the results of two scenario. We observe that, the absolute bias and MSE of the estimators of $h_j(\cdot)$ for $j=1,2$ decrease with increase in sample size.


\begin{table}[h]\label{ph-h1andh2-1}
\caption{Absolute bias and MSE of the estimators of $h_1(.)$ and $h_2(.)$ under PH model at various time points for $\beta_1=2, \beta_2=2$}
\label{t:one}
\begin{small}
\centering
\begin{tabular}{crrrrrrrrr}

	\hline {Censoring(\%)}  & $ n$
	&Bias &MSE &Bias &MSE&Bias&MSE\\
	\hline
	&&\multicolumn{2}{c}{$t_1=1$}&\multicolumn{2}{c}{$t_2=2$}&
\multicolumn{2}{c}{$t_3=3$}\\
		\hline
		&&\multicolumn{6}{c}{$h_1(.)$}\\	
	\hline	
	\multirow{3}{*}{20} & 100 & 0.1334&	0.0229&	0.1482&	0.0279	&0.1046&	0.0172 \\
	
	& 200 & 0.0119&	0.0103&	0.0552&	0.0152&	0.0378&	0.0142 \\
	& 500 & 0.0032&	0.0021&	0.0138&	0.0025&	0.0056&	0.0024  \\
	\hline
	\multirow{3}{*}{40} & 100& 0.0661&	0.0110&	0.1008&	0.0230&	0.0840&	0.0206
\\
	& 200 & 0.0266&	0.0096&	0.0502&	0.0086&	0.0601&	0.0101 	\\
	& 500 &0.0024&	0.0028&	0.0049&	0.0025&	0.0012&	0.0026 	\\

		\hline
		&&\multicolumn{6}{c}{$h_2(.)$}\\	
		\hline
\multirow{3}{*}{20} & 100 &0.0721&	0.0206&	0.0858&	0.0194&	0.1062&	0.0271\\

& 200 & 0.0553&	0.0101	&0.0643&	0.0125&	0.0944&	0.0199 \\
& 500 &0.0034&	0.0028&	0.0452&	0.0106&	0.0885&	0.0112 \\
\hline
\multirow{3}{*}{40} & 100& 0.2594&	0.0879&	0.2479	&0.0853	&0.2187&0.0727
\\
& 200 &0.2594&	0.0879&	0.2479&	0.0853&	0.2187&	0.0727 \\
& 500 &0.0260&	0.0045&	0.0279&	0.0053&	0.0166&	0.0050 \\
\hline

\end{tabular}
\end{small}
\end{table}%


\begin{table}[h]\label{ph-h1andh2-2}
\caption{ Absolute bias and MSE of the estimators of $h_1(.)$ and $h_2(.)$ under PH model at various time points for $\beta_1=3, \beta_2=2$}
\label{t:one}
\begin{small}
\centering
\begin{tabular}{crrrrrrr}

	\hline {Censoring(\%)}  & $ n$
	&Bias &MSE &Bias &MSE&Bias&MSE\\
	\hline
	&&\multicolumn{2}{c}{$t_1=1$}&\multicolumn{2}{c}{$t_2=2$}&
\multicolumn{2}{c}{$t_3=3$}\\
		\hline
		&&\multicolumn{6}{c}{$h_1(.)$}\\	
	\hline	
	\multirow{3}{*}{20} & 100 &  0.0590	&0.0148&0.0765&	0.0170&	0.0811&	0.0190
\\

	& 200 & 0.0324&	0.0071&	0.0172&	0.0092&	0.0364&	0.0101 \\
	& 500 & 0.0124&	0.0063&	0.0028&	0.0083&	0.0107&	0.0088  \\
	\hline
	\multirow{3}{*}{40} & 100& 0.1293&	0.0241&	0.1785&	0.0479&	0.1922&	0.0565

\\
	& 200 & 0.0634&	0.0192&	0.1747&	0.0394&	0.1767&	0.0405 	\\
	& 500 & 0.0148&	0.0032&	0.0287&	0.0043&	0.0350&	0.0048 	\\

		\hline
		&&\multicolumn{6}{c}{$h_2(.)$}\\	

\hline
\multirow{3}{*}{20} & 100 &0.1358&	0.0256&	0.0918&	0.0171&	0.1468&	0.0303
\\

& 200 & 0.0278&	0.0152&	0.0659&	0.0168&	0.0529&	0.0185 \\
& 500 & 0.0266&	0.0035&	0.0090&	0.0075&	0.0251&	0.0061 \\
\hline
\multirow{3}{*}{40} & 100& 0.1718&	0.0405&	0.2640&	0.0771&	0.2714&	0.1007
\\
& 200 &0.1357&	0.0396&	0.2552&	0.0715	&0.2494	&0.0747 \\
& 500 & 0.0106&	0.0040&	0.0124&	0.0047&	0.0052	&0.0048 \\
\hline
\end{tabular}
\end{small}
\end{table}%


\begin{table}[h]\label{po-h1andh2-1}
\caption{ Absolute bias and MSE of the estimators of $h_1(.)$ and $h_2(.)$ under PO model at various time points for $\beta_1=2, \beta_2=3$}
\label{t:one}
\begin{small}
\centering
\begin{tabular}{crrrrrrr}

	\hline {Censoring(\%)}  & $ n$
	&Bias &MSE &Bias &MSE&Bias&MSE\\
	\hline
	&&\multicolumn{2}{c}{$t_1=1$}&\multicolumn{2}{c}{$t_2=2$}&
\multicolumn{2}{c}{$t_3=3$}\\
		\hline
		&&\multicolumn{6}{c}{$h_1(.)$}\\	
	\hline	
	\multirow{3}{*}{20} & 100 &  0.0611	&0.0140	&0.0454	&0.0124&0.0367&	0.0140
\\
	
	& 200 & 0.0249&	0.0051&	0.0264&	0.0065&	0.0325&	0.0061 \\
	& 500 & 0.0144	&0.0026	&0.0034	&0.0044&0.0027&	0.0038  \\
	\hline
	\multirow{3}{*}{40} & 100& 0.0253&	0.0117&	0.0562&	0.0163&	0.0480&	0.0141
\\
	& 200 & 0.0176	&0.0057	&0.0392&0.0064&	0.0186&	0.0090 	\\
	& 500 &	0.0079&	0.0022&	0.0067&	0.0040&	0.0172&	0.0029 \\

		\hline
		&&\multicolumn{6}{c}{$h_2(.)$}\\	
		\hline

\multirow{3}{*}{20} & 100 & 0.1257&	0.0357&	0.1559&	0.0358&	0.1709&	0.0413
\\

& 200 & 0.1234&	0.0249&	0.0595&	0.0242&	0.0717&	0.0256  \\
& 500 & 0.0673&	0.0083	&0.0222&	0.0080&	0.0278&	0.0098 \\
\hline
\multirow{3}{*}{40} & 100& 0.1626&	0.0396&	0.0977&	0.0379	&0.1011	&0.0446
\\
& 200 &0.1081&	0.0395&	0.0703&	0.0249&	0.0728&	0.0214 \\
& 500 & 0.0093&	0.0050&	0.0123&	0.0060&	0.0009&	0.0061 \\
\hline

\end{tabular}
\end{small}
\end{table}%


\begin{table}[h]\label{po-h1andh2-2}
\caption{ Absolute bias and MSE of the estimators of $h_1(.)$ and $h_2(.)$ under PO model at various time points for $\beta_1=3, \beta_2=3$}
\label{t:one}
\begin{small}
\centering
\begin{tabular}{crrrrrrr}

	\hline {Censoring(\%)}  & $ n$
	&Bias &MSE &Bias &MS2&Bias&MSE\\
	\hline
	&&\multicolumn{2}{c}{$t_1=1$}&\multicolumn{2}{c}{$t_2=2$}&
\multicolumn{2}{c}{$t_3=3$}\\
		\hline
		&&\multicolumn{6}{c}{$h_1(.)$}\\	
	\hline	
	\multirow{3}{*}{20} & 100 &  0.0809& 	0.0206& 	0.1048& 	0.0184& 	0.0861& 	0.0185 \\
	& 200 &  0.0720& 	0.0080& 0.0354& 0.0143& 0.0215& 0.0109 \\
	& 500 &  0.0138	& 0.0069& 0.0100& 0.0082& 0.0100& 0.0086 \\
	\hline
	\multirow{3}{*}{40} & 100& 0.1567& 	0.0400& 0.1172& 0.0321& 0.1007& 	0.0295 \\
	& 200 &  0.0136& 0.0074& 0.0398& 0.0104	& 0.0639&  0.0109 	\\
	& 500 &	 0.0099 & 0.0031& 0.0384& 0.0050& 0.0414& 0.0078 \\

		\hline
		&&\multicolumn{6}{c}{$h_2(.)$}\\	
\hline
\multirow{3}{*}{20} & 100 & 0.1434&	0.0415&	0.2016&	0.0652	&0.2140	&0.0716
\\
& 200 & 0.1231&	0.0249&	0.0884&	0.0193&	0.0767&	0.0179 \\
& 500 &  0.0066	&0.0037	&0.0031&0.0044&	0.0101&0.0047 \\
\hline
\multirow{3}{*}{40} & 100& 0.0696&0.0325&	0.1923&	0.0693&	0.1645&	0.0448
\\
& 200 & 0.0348&0.0144&	0.1453&	0.0366&	0.1052&	0.0434 \\
& 500 & 0.0238&	0.0055&	0.0276&	0.0068&	0.0378&	0.0077  \\
\hline

\end{tabular}
\end{small}
\end{table}%

\begin{table}[h]\label{ph-cp and aw}
\caption{ Coverage probability and average width of the confidence interval of the estimates of $\beta_1$ and $\beta_2$ under PH model}
\label{t:one}
\begin{small}
\centering
\begin{tabular}{crrrrrrrrr}

	\hline {Censoring(\%)}  & $ n$
	&CP&AW&CP&AW&CP&AW&CP&AW\\
		\hline
		&&\multicolumn{2}{c}{$\beta_1=2$}&\multicolumn{2}{c}{$\beta_2=2$}&
\multicolumn{2}{c}{$\beta_1=2$}&\multicolumn{2}{c}{$\beta_2=3$}\\	
	\hline	
	\multirow{3}{*}{20} & 100 & 0.9507&	0.5628&	0.9529&	0.6710&0.9529&	0.5606&	0.9524&	0.7391 \\
	& 200 & 0.9504&	0.3918&	0.9518&	0.4564 &0.9528&	0.3886&	0.9506&	0.5125 \\
	& 500 & 0.9497&	0.2443&	0.9487&	0.2882& 0.9499&	0.2457&	0.9495&	0.3258 	\\
	\hline
	\multirow{3}{*}{40} & 100&  0.9529&	0.5751&	0.9539&	0.7796 &0.9530	&0.5763&	0.9547&	0.8544 \\
	& 200 & 0.9526&	0.4012&	0.9522&	0.5385 &0.9487&	0.4055&	0.9505&	0.5911 	\\
	& 500 &0.9509&	0.2547&	0.9512&	0.3440& 0.9494&	0.2518	&0.9498&0.3703 	\\

		\hline
		&&\multicolumn{2}{c}{$\beta_1=3$}&\multicolumn{2}{c}{$\beta_2=2$}&
\multicolumn{2}{c}{$\beta_1=3$}&\multicolumn{2}{c}{$\beta_2=3$}\\
		\hline
\multirow{3}{*}{20} & 100 &0.9545&	0.6212&	0.9533&	0.6563& 0.9514&	0.6219&	0.9538&	0.7406 \\
& 200 & 0.9501&	0.4393&	0.9517&	0.4654 &0.9521&	0.4401&	0.9512&	0.5147 \\
& 500 & 0.9526&	0.2789&	0.9491&	0.2902& 0.9517	&0.2736&0.9521&	0.3222  \\
\hline
\multirow{3}{*}{40} & 100& 0.9532&	0.6488&	0.9555&	0.8069 &0.9527&	0.6522&	0.9535&	0.8745 \\
& 200 &0.9507&	0.4519&	0.9518&	0.5488&0.9516&	0.4593&	0.9537&	0.6000 \\
& 500 &0.9499&	0.2864&	0.2864&	0.3451& 0.9518&	0.2888&	0.9509	&0.3721 \\
\hline
\end{tabular}
\end{small}
\end{table}%


\begin{table}[h]\label{po-cp and aw}
\caption{ Coverage probability and average width of the confidence interval of the estimates of $\beta_1$ and $\beta_2$ under PO model}
\label{t:one}
\begin{small}
\centering
\begin{tabular}{crrrrrrrrr}

	\hline {Censoring(\%)}  & $ n$
	&CP&AW&CP&AW&CP&AW&CP&AW\\
		\hline
		&&\multicolumn{2}{c}{$\beta_1=2$}&\multicolumn{2}{c}{$\beta_2=2$}&
\multicolumn{2}{c}{$\beta_1=2$}&\multicolumn{2}{c}{$\beta_2=3$}\\	
	\hline	
	\multirow{3}{*}{20} & 100 & 0.9523	&0.5546	&0.9523&	0.6599 &0.9520&	0.5567&	0.9542	&0.7385  \\
	& 200 & 0.9507&	0.3917&	0.9521&	0.4621& 0.9515&	0.3863&	0.9520&	0.5182  \\
	& 500 & 0.9499	&0.2460	&0.9491	&0.9491& 0.9511	&0.2474&0.9506&	0.3226 \\
	\hline
	\multirow{3}{*}{40} & 100&  0.9491&	0.5766	&0.9561	&0.7904 &0.9527	&0.5800&	0.9539&	0.8649  \\
	& 200 &  0.9493	&0.4063	&0.9519	&0.5472& 0.9514	&0.4096	&0.9524	&0.5888 \\
	& 500 & 0.9508&	0.2540&	0.9485&	0.3394& 0.9512&	0.2540&	0.9508&	0.3657 \\

		\hline
		&&\multicolumn{2}{c}{$\beta_1=3$}&\multicolumn{2}{c}{$\beta_2=2$}&
\multicolumn{2}{c}{$\beta_1=3$}&\multicolumn{2}{c}{$\beta_2=3$}\\
		\hline
\multirow{3}{*}{20} & 100 &  0.9541	&0.6232	&0.9528	&0.6604 &0.9530	&0.6318&	0.9524&	0.7407 \\
& 200 & 0.9538 &0.4372&	0.9509	&0.4625 &0.9474&0.4353&	0.9525&	0.5133  \\
& 500 & 0.9508&	0.2742	&0.9496	&0.2888 &0.9505&0.2788&	0.9484&	0.3195 \\
\hline
\multirow{3}{*}{40} & 100&0.9522&0.6453&0.9537&	0.7964 &0.9491	&0.6463&	0.9537	&0.8690  \\
& 200 & 0.9521&	0.4516	&0.9518&0.5503& 0.9508	&0.4500	&0.9530	&0.5926 \\
& 500 & 0.9504&	0.2865&	0.9507&	0.3388 &0.9508&	0.2844	&0.9511	&0.3713 \\
\hline
\end{tabular}
\end{small}
\end{table}%
Next, we obtain the coverage probability and average width of the confidence interval of the regression parameters for different parameter combinations considered above. The results for PH and PO models are given in Tables 7 and 8, respectively. We observe that, the coverage probability approaches 0.95 and the average width of the interval decreases as $n$ increases for both PH and PO models.

Finally, we compare theoretical asymptotic variance and Monte Carlo variance of the estimators of $\beta_1$ and $\beta_2$ for different parameter combinations under PH model. In PH model, we obtain the variance-covariance matrix specified in  (\ref{var}) and (\ref{limvar}) as $$\Sigma_{0k}=\Sigma_{1k}=Var\left(\int_{0}^{\infty}(Z-\mu_k(t))dM_k(t)\right), \, k=1,2,$$where $\mu_k(t)=E(Z|\tilde T,\delta=1)$. Hence $\Sigma_k=\Sigma_{0k}^{-1}$. For more details see Chen, Jin and Zhang (2002). Comparison results are given in Table 9.
We observe that the estimators of theoretical variance (EAV) and the estimators of Monte Carlo variance (MCV), agree each other.


\begin{table}[h]\label{ph-variances}
\caption{ Theoretical asymptotic variance and Monte Carlo variance of  the estimates of $\beta_1$ and $\beta_2$ under PH model}
\label{t:one}
\begin{small}
\centering
\begin{tabular}{crrrrrrrrr}

	\hline {Censoring(\%)}  & $ n$
	&EAV&MCV&EAV&MCV&EAV&MCV&EAV&MCV\\
		\hline
		&&\multicolumn{2}{c}{$\beta_1=2$}&\multicolumn{2}{c}{$\beta_2=2$}&
\multicolumn{2}{c}{$\beta_1=2$}&\multicolumn{2}{c}{$\beta_2=3$}\\	
	\hline	
	\multirow{3}{*}{20} & 100 & 0.0200&	0.0206&	0.0221&	0.0273& 0.0200&	0.0204&	0.0200&	0.0356\\
	& 200 & 0.0100&	0.0100&	0.0100&	0.0136  &0.0100&0.0098&	0.0100&0.0171  \\
	& 500 &  0.0040	&0.0039&0.0040&	0.0054&0.0040&	0.0039&	0.0040&	0.0069 	\\
	\hline
	\multirow{3}{*}{40} & 100&0.0200&	0.0215&	0.0201&	0.0395& 0.0200&	0.0216&	0.0201&	0.0475  \\
	& 200 &0.0100&	0.0105&	0.0100&	0.0189& 0.0100&	0.0107&	0.0100&0.0227 	\\
	& 500 & 0.0040	 & 0.0042 & 0.0040 & 0.0077 & 0.0040 & 0.0041 & 0.0040 & 	0.0089 	\\

		\hline
		&&\multicolumn{2}{c}{$\beta_1=3$}&\multicolumn{2}{c}{$\beta_2=2$}&
\multicolumn{2}{c}{$\beta_1=3$}&\multicolumn{2}{c}{$\beta_2=3$}\\
		\hline
\multirow{3}{*}{20} & 100 & 0.0200&	0.0251&	0.0200&	0.0280& 0.0200&	0.0252&	0.0200&	0.0357 \\
& 200 & 0.0100&	0.0126&	0.0100&	0.0141 & 0.0100	&0.0126&0.0100	&0.0172  \\
& 500 & 0.0040&	0.0051&	0.0040&	0.0055 & 0.0040	&0.0049	&0.0040&0.0068   \\
\hline
\multirow{3}{*}{40} & 100& 0.0200&	0.0274&	0.0201&	0.0424&  0.0200	&0.0277&	0.0201	&0.0498   \\
& 200 &0.0100&0.0133&0.0100	&0.0196 & 0.0100&0.0137	&0.0100	&0.0234   \\
& 500 & 0.3451&	0.0053&	0.0040&	0.0077&  0.0040&0.0054&	0.0040&	0.0090 \\
\hline
\end{tabular}
\end{small}
\end{table}%

\section{An example}
In this section, we apply our proposed method to the data obtained from a clinical trial on HIV infection and AIDS of 329 homosexual men from Amsterdam. The data is available in R package ‘mstate’ and is exclusively studied by Geskus (2015). During the course from HIV infection from non-syncytium-inducing (NSI) phenotype to death, intermediate events may occur that have an impact on subsequent disease progression. One such event is a switch of the HIV virus to the syncytium inducing (SI) phenotype and the other is the progression to AIDS.  Since neither of the events is final and changes the probability of relapse/ survival of a patient, it is of interest to know whether AIDS or SI is the first event to occur.  As mentioned in the introductory section,  the data may contain cured individuals.  Accordingly, we use this data to illustrate our newly proposed method. We use the estimating equation (\ref{betahat}) for finding the estimators of regression coefficients. We estimate cumulative incidence functions for patients with AIDS and SI using proportional hazard (PH) model and proportional odds (PO) models.

The presence or absence of the CCR5-$\Delta$32 deletion in one or both chromosomes affect the progression to AIDS or SI significantly, it is considered as one covariate in this study. We also consider the age of the patient at the time of infection of HIV as another covariate. Individuals without the deletion of one of the chromosomes are referred to as WW (wild type allele on both chromosomes) and those who have the deletion are considered as having a mutation and referred to as WM (mutant allele on one chromosome). We removed 5 records due to missingness of covariate value CCR5-$\Delta$32. Out of 324 patients, the first event to occur was AIDS for 113 patients and SI was occurred as the first event to 107 patients. The remaining 104 patients are observed to be event free in the study period.  For 259 patients the covariate value is `WW’ and for the remaining 65 patients it is `WM’.
An interesting characteristic of the data is that  time to occurrence of the event varies from 0.112 years (approximately 41 days) to 13.936 years (approximately 5090 days).  We also see that after 12.4 years (4526 days) there are only 3 event times observed one due to AIDS and two with the cause SI, whereas around 14.5\% of the total lifetimes are larger than 12.4 years. Presence of  large number of right censored observations  indicates the possible presence of cured individuals in the population. This fact motivates us to analyse the data using the proposed method. To ensure the presence of cured individuals in the population, we plot Kaplan Meier curve of the complete data in Figure 1. From Figure 1, we can see that the minimum survival probability is around 0.2, and then the curve experiences a sudden fall. This is due to the fact that the largest time is an observed lifetime. To differentiate between the behavior of causes AIDS and SI over time, we also plot the baseline cumulative incidence functions in Figure 2.
\begin{figure}[h]
	\centering
	\includegraphics[height=10cm,width=120mm]{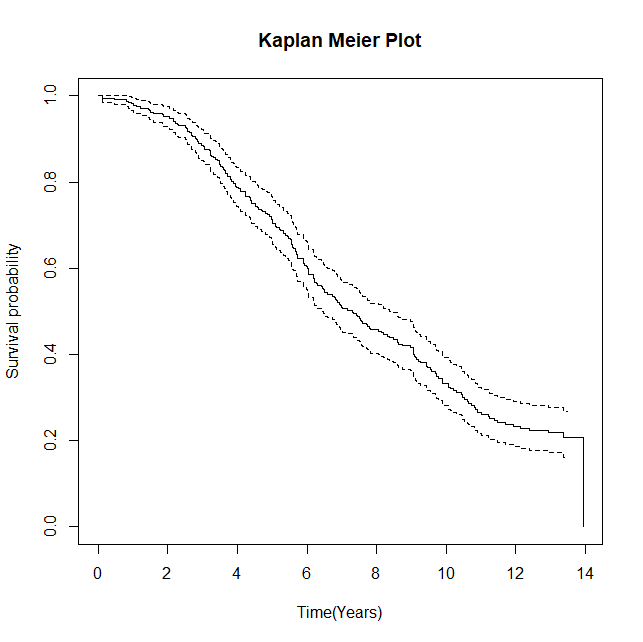}
	\caption{Kaplan Meier curve of the complete data}
\end{figure}

\begin{figure}[h]
	\centering
	\includegraphics[height=10cm,width=120mm]{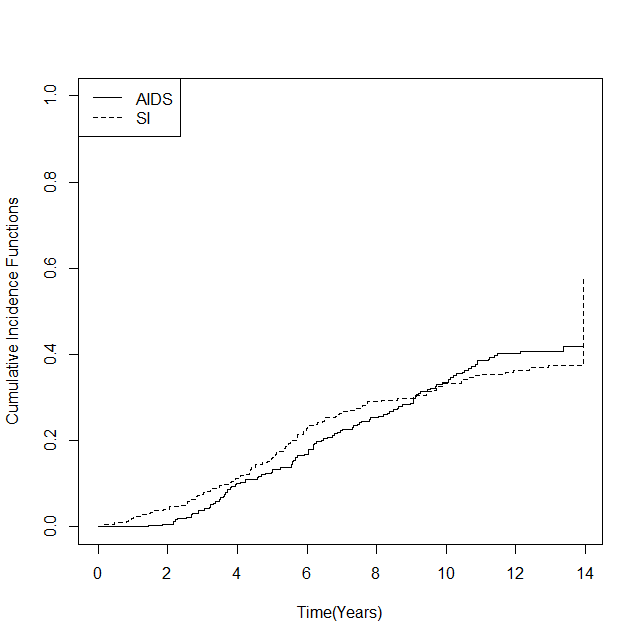}
	\caption{ Baseline cumulative incidence functions due to AIDS and SI}
\end{figure}

The regression coefficients are estimated using the estimating equation (\ref{betahat}) under  PH  and PO model assumption. We use the R package ` `TransModel' to estimate the regression coefficients.  The estimated  regression coefficients  with corresponding standard errors are reported  in Table 10.
\begin{table}
\label{t:five}
		\centering
	\caption{ Estimates of the regression parameters with corresponding standard error}
	\begin{small}
				\begin{center}
			\begin{tabular}{ cccccc }
				\hline
				Cause & Model  & Covariate & Coefficient & SE & P-value\\
				\hline
				AIDS & PH & CCR5 &-1.2001  & 0.3081 & 0.0005\\
					 &  & Age &0.0185  & 0.0141 & 0.1908\\
				&PO& CCR5& -1.5439& 0.3619& 0.0005 \\
				&& Age&0.0241&0.0184& 0.1905 \\
				\hline
				SI & PH & CCR5 &-0.2284  & 0.2390 & 0.3393\\
					 &  & Age &0.0156  & 0.0140 & 0.2660\\
				&PO& CCR5&-0.2942&0.2991&  0.3255\\
				&& Age&0.0193&0.0173&  0.2646\\
				\hline
			\end{tabular}
		\end{center}
	\end{small}
\end{table}
To understand how the covariate values affect the lifetimes due to different  causes, we plot  the cumulative incidence functions due to AIDS and SI separately. To plot the graphs, we categorize the covariate age as less than or equal to 34.5 years (median age) and grater  than 34.5 years. We plot the graphs of cumulative incidence functions of AIDS patients, for various categories according to the covariate values in Figure 3 and the same for SI patients in Figure 4.
Cure fraction is estimated using Eq. (\ref{cureprobest}). The value of cure fraction  obtained using PH $0.1655$. This estimated value of cure fraction also supports the claim that the data consist of cured individuals as evident from the Kaplan Meier curve. As mentioned earlier,  here it is not required to  model the $\pi(Z)$'s separately to find the overall cure fraction.
To compare our estimator of the cure rate with the existing methods, we fit a logistic regression  model (Farewell, 1982) for the  data by ignoring causes and then estimate the cure rate. We made the lifetime as binary variable by taking 12.4 years as cutoff.   The estimate of the cure rate is obtained as $0.2961$. This might be an over estimate of the cure rate as evident from the discussion above.

 \begin{figure}[h]
	\centering
	\includegraphics[height=10cm,width=120mm]{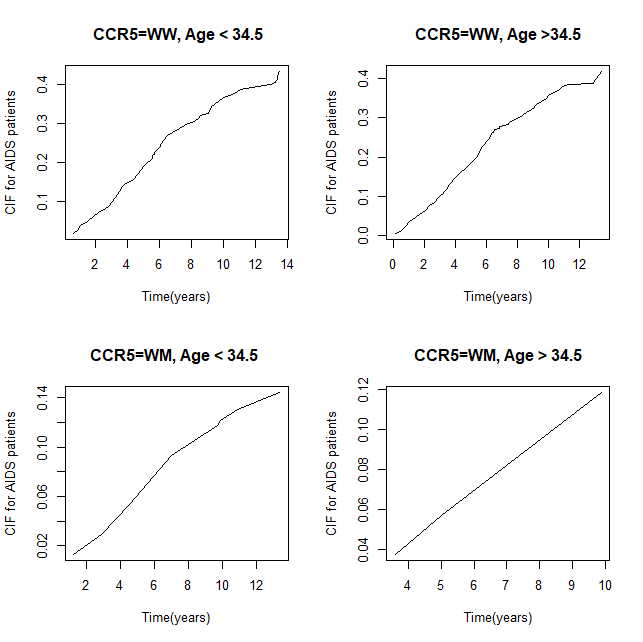}
	\caption{ CIF of patients with AIDS under PH model in different categories }
 \end{figure}

\begin{figure}[h]
	\centering
	\includegraphics[height=10cm,width=120mm]{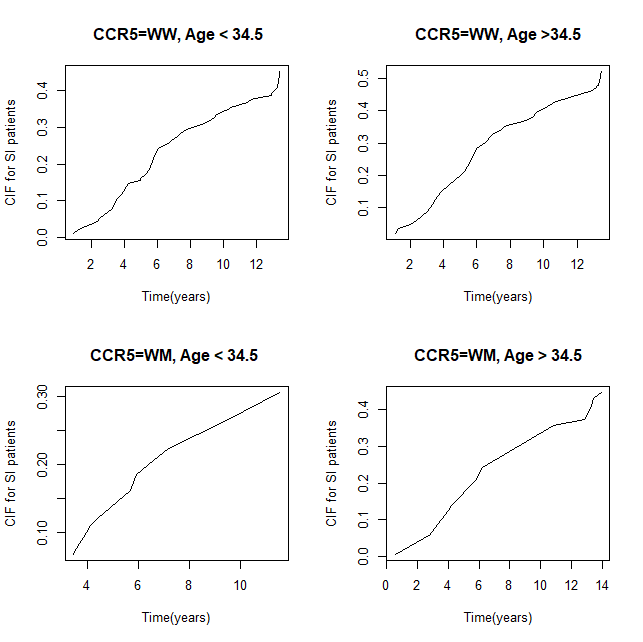}
	\caption{ CIF of patients with SI under PH model in different categories }
 \end{figure}

\section{Concluding Remarks}
Cure rate models become popular due to its applications in various  fields.  In the present study, we proposed a  semi-parametric model for the analysis of competing risks data with cure fraction. We developed estimators of cumulative incidence  functions using linear transformation models. Unlike the traditional approach, it does not require  to find the estimator of the cure fraction to compute the overall survival function, which leads to a simple procedure. The regression parameters are estimated using martingale based estimating equations. The asymptotic distribution of the estimators was shown to be Gaussian.  The finite sample  performance of the estimators of $\beta_k$ and and $h_k(\cdot)$ at different time points, are evaluated in terms of bias and MSE, through  a Monte Carlo simulation study. We also find the coverage probability and the  average width of the regression estimators.  A real data  obtained from a study of HIV and AIDS infection was analysed using the proposed method. The newly developed method has a good impact  due to the model flexibility and computational advantages.

In the present study, we considered right censored data.  Different types of censoring schemes such as current status censoring, double censoring and interval censoring are common in survival studies. The proposed method can be extended to these set up  by suitably constructing martingale based estimating equations. The works in this direction will be reported elsewhere. The model selection is an important concern in  the analysis of competing risks data with cure fraction when the cumulative incidence functions are specified through  linear transformation. A criterion for model selection can be developed using martingale based residuals. This problem has to be addressed separately.
 \section*{Acknowledgements}Sudheesh K. K. is thankful to SERB, DST, Govt. of India for financial support received under MATRICS grant. Sreedevi E. P. and Sankaran P. G. would like to thank Kerala State Council for Science Technology  and Environment, Kerala, India for the financial support provided to carry out this research work. We extend our thanks to anonymous reviewers for their constructive suggestions.

\vspace{-0.2in}
\appendix
\section{}\vspace{-0.2in}
\subsection{Proof of Theorem 1:} The estimators of $F_k(t|Z)$ are obtained from equation (\ref{LTMCOMPSUR}) by replacing ${\beta}_k$ and ${h}_k(t)$ with its estimators. Then we obtain
$$\widehat F_k(t|Z)=g_k^{-1}(\widehat h_k(t) + Z'\widehat\beta_k),\,\,k=1,\ldots,K.$$
 In cure rate model, we know that $\lim_{t\rightarrow \infty}S(t|Z)=a_0>0$.  Hence in view of  the relation specified in equation (\ref{mixsurnew}), we have $\liminf\big\{1-\sum_{k=1}^{K}\widehat F_k(t|Z)\big\}>0$. The estimators  $\widehat{\beta}_k$ and $\widehat{h}_k(t)$ are obtained by solving the equations (\ref{betahat}), (\ref{eq5.1}) and (\ref{eq5.2}). Now using the assumption $h_k(0)=-\infty$ and that $h_k(t)$ is non-decreasing function, from these three equations it follows that, $\widehat{\Lambda}_{\epsilon_{k}}(\tau)<\infty$, where $\tau$ is defined in Assumption $D3$.  Therefore, the consistency of the estimators $\widehat{\beta}_k$ and $\widehat{\Lambda}_{\epsilon_{k}}(t)$ follows from Theorem 1 of Mao and Lin (2017). By Step A1 of Chen et al. (2002, p.665), we have consistency of $\widehat{h}_k(t)$, which completes the proof of the Theorem 1.
\vspace{-0.2in}
\subsection{Proof of Theorem 2:}  We use martingale central limit theorem to establish  the asymptotic normality  of $\widehat{\beta}_k$, $k=1,\ldots,K$.

First we represents   $U(\widehat{\beta}_k,\widehat{h}_k(t))$ as
$$U\big(\widehat{\beta}_k,\widehat{h}_k(t)\big)=U\big({\beta}_k,\widehat{h}_k(t)\big)+
(\widehat{\beta}_k-\beta_k)\frac{\partial}{\partial \beta_k}U\big({\beta}_k,\widehat{h}_k(t)\big).$$
Since $\widehat{\beta}_k$ and $\widehat{h}_k(t)$ are the solutions of equations (\ref{betahat}) and (\ref{LTMO}) we have    $U(\widehat{\beta}_k,\widehat{h}_k(t))=0$. Now we can write
 \begin{equation}\label{beta}
\sqrt{n}(\widehat{\beta}_k-\beta_k)=(\frac{1}{n}\frac{\partial}{\partial \beta_k}U\big({\beta}_k,\widehat{h}_k(t)\big))^{-1}
\frac{1}{\sqrt{n}}U\big({\beta}_k,\widehat{h}_k(t)\big).
\end{equation}To establish the asymptotic normality of $\sqrt{n}(\widehat{\beta}_k-\beta_k)$, first we show that, as $n\rightarrow\infty$,    $(\frac{1}{n}\frac{\partial}{\partial \beta_k}U\big({\beta}_k,\widehat{h}_k(t)\big) )^{-1}$ converges in probability to $\Sigma_{1k}$ given in Eq. (\ref{limvar})  and $ \frac{1}{\sqrt{n}}U\big({\beta}_k,\widehat{h}_k(t)\big)$ converges in distribution to Gaussian.  Then  the asymptotic normality of $\sqrt{n}(\widehat{\beta}_k-\beta_k)$ follows from  equation (\ref{beta}) by  applying  Slutsky's theorem.

Using similar argument of Step A3 of Chen et al. (2002), we have
\begin{equation}\label{lln}
  \frac{1}{n}\frac{\partial}{\partial \beta_k}U\big({\beta}_k,\widehat{h}_k(t)\big) =\Sigma_{1k}+o_p(1),
\end{equation}
where $\Sigma_{1k}$ is specified in Eq. (\ref{limvar}).
In view of the martingale representation given in  (\ref{martingale}), using the Step A4 of Chen et al. (2002), we have
$$ \frac{1}{\sqrt{n}}U\big({\beta}_k,\widehat{h}_k(t)\big)=\frac{1}{\sqrt n}\sum_{i=1}^{n}\int_{0}^{\tau}({Z}-\mu_k(t))dM_{ik}(t)+o_p(1).$$
Using martingale central limit theorem, as $n\rightarrow\infty$, $ \frac{1}{\sqrt{n}}U\big({\beta}_k,\widehat{h}_k(t)\big)$ converges in distribution  to multivariate normal with mean vector  zero and variance-covariance matrix  $\Sigma_{0k}$, where  $\Sigma_{0k}$ is the limit of the predictable variation process. Using the martingale representation given in (\ref{martingale}), the predictable variation process is given by
$$\Sigma_{0k}^{*}=\frac{1}{ n}\sum_{i=1}^{n}\int_{0}^{\tau}({Z}-\mu_k(t))'({Z}-\mu_k(t))Y(t)d\Lambda_{\epsilon_{ik}}(Z'\beta_{k}+{h}_{k}(t)).$$
By law of large numbers, as $n\rightarrow\infty$, $\Sigma_{0k}^{*}$ converges in probability to $\Sigma_{0k}$, where $\Sigma_{0k}$ is  specified in Eq. (\ref{var}). In view of the representations  (\ref{beta}) and (\ref{lln})  by applying  Slutsky's theorem,  as $n\rightarrow\infty$,  $\sqrt{n}(\widehat{\beta}_k-\beta_k)$ converges in distribution to multivariate normal with    zero mean vector  and variance-covariance matrix  $\Sigma_k=\Sigma_{1k}^{-1}\Sigma_{0k}(\Sigma_{1k}^{-1})'$. This completes the proof of the Theorem 2.

\end{document}